
\def\picture #1 by #2 (#3){$
\vcenter to #2{
\hrule width #1 height 0pt depth 0pt
\vfill
\special{picture #3}}$}

\documentstyle{amsppt}
\nologo

\centerline{\tenbf THE PRINCIPLE DESCRIBING POSSIBLE COMBINATIONS
OF}
\vskip 5.5pt
\topmatter
\title Singularities in deformations of a fixed singularity\endtitle
\author
Tohsuke Urabe
\endauthor\affil
Department of Mathematics\\
Tokyo Metropolitan University\\
Minami-Ohsawa 1-1,\ Hachioji-shi\\
Tokyo 192-03 \ Japan\\
(E-mail: urabe\@ math.metro-u.ac.jp)\endaffil
\endtopmatter

\document
This is a review of my recent work.

What kinds of combinations of singularities can appear in small
deformation fibers of a fixed singularity?
We consider this problem for hypersurface singularities on complex
analytic spaces of dimension 2.

Recall first that a connected Dynkin graph of type $A$, $D$  or
$E$ corresponds to a surface singularity called a {\it rational
double point}. (Durfee \cite{3}.)

When the fixed singularity is a rational double point
corresponding to a Dynkin graph $\Gamma_0$, the answer to our
problem is well-known.
\roster
\item  Any small deformation fiber has only rational double points
as singularities.
\item A combination of rational double points corresponding
to a Dynkin graph $\Gamma$ appears on a small deformation fiber
if and only if $\Gamma$ is a subgraph of $\Gamma_0$.  (The type of
each component of
$\Gamma$  corresponds to the type of the singularity on  a small
deformation fiber $Y$ and the number of components of each type
corresponds to the number of singularities of each type on  $Y$.)
\endroster
This follows from the description of the semi-universal deformation
family of a rational double point due to Grothendieck and
Brieskorn. In this case the deformation family can be identified
with a subspace of the corresponding simple Lie algebra, and the
monodromy covering of the parameter space of the family can be
identified with the Cartan subalgebra of the corresponding type.
(Slodowy \cite{6}.) Thus the above fact follows.

We would like to consider more complicated singularity in the higher
hierarchies in Arnold's classification list. (Arnold \cite{1}. There
we can find the defining polynomials of singularities
considered below.)

Let  $\Xi$   be a class of surface singularities.  By  $PC(\Xi)$ we
denote the set of Dynkin graphs $\Gamma$ with several components
such that there exists a small deformation fiber  $Y$  of a
singularity belonging to  $\Xi$  satisfying the following
conditions:

\roster
\item  $Y$  has only rational double points as singularities.
\item The combination of rational double points on  $Y$  corresponds
exactly to  $\Gamma$.
\endroster

When $\Xi$ is one of three kinds $P_8$, $X_9$, $J_{10}$ of simple
elliptic singularities, the answer to the above problem is the
following: If a small deformation fiber  has a singularity
which is not a rational double point, then the singularity
is unique and it is a simple elliptic singularity of the
same type $\Xi$. Besides, $\Gamma$ belongs to
$PC(\Xi)$ if and only if
$\Gamma$ can be made by elementary transformations repeated twice
from the corresponding basic graph $\Gamma_0$. The basic graph
$\Gamma_0$ is the Dynkin graph of type $E_6$ if $\Xi=P_8$, of type
$E_7$ if
$\Xi=X_9$, or of type $E_8$ if $\Xi=J_{10}$. (Looijenga \cite{5},
Urabe \cite{7}.)

An {\it elementary transformation} is an operation by which we can
make a new Dynkin graph from a given Dynkin graph.  We give the
definition below.

\proclaim{ Definition}\rm
 (An elementary transformation)\quad The following
procedure  is called an elementary transformation of a Dynkin
graph.

\roster
\item Replace each connected component by the corresponding extended
Dynkin graph.

\item Choose in arbitrary manner at least one vertex from each
component (of the extended Dynkin graph) and then remove these
vertices together with edges issuing from them.
\endroster\endproclaim

We can find the definition of the extended Dynkin graph in any book
on Lie algebras. (Bourbaki \cite{2}.) They can be made by adding one
vertex and one or two edges to each connected component of the
Dynkin graph. The position of the added vertex and edges depends on
the type of the component.

\proclaim{Example}\rm
We start from $E_7$. Removing a vertex of the extended Dynkin graph
of type $E_7$ as in the following figure, we can
make the graph $D_6+A_1$:

\centerline{\picture 101mm by 25mm (E7 scaled 1000)}

Applying another elementary transformation to $D_6+A_1$, we can
make the graph $D_4+3A_1$.

\centerline{\picture 107mm by 35mm (D6A1 scaled 1000)}

Since $E_7$ is the basic graph for
$\Xi=X_9$, we have $D_4+3A_1\in PC(X_9)$
\endproclaim

In the case where $\Xi$ is the class of
cusp singularities
$x^p+y^q+z^r+\lambda xyz=0$ of type $(p,\,q,\,r)$, (The indices $p$,
$q$ and $r$ are positive integers with
$1/p+1/q+1/r<1$. We may assume $p\le q\le r$. The parameter
$\lambda$ is a non-zero constant.) Loijenga gave the answer. We call
the graph below the Gabri\'{e}lov graph of type
$(p,\,q,\,r)$. The numbers of vertices in the three arms are $p$,
$q$ and $r$ respectively including the common central one. A subgraph
is called a Dynkin subgraph, if all of its components are Dynkin
graphs of type
$A$,
$D$ or $E$.

\centerline{\picture 121mm by 40mm (Gabrielov scaled 1000)}

\vskip 10pt

According to Loijenga, $\Gamma\in PC(\Xi)$ if and only if $\Gamma$
can be made by one elementary transformation from a Dynkin subgraph
of the  Gabri\'{e}lov graph of type $(p,\,q,\,r)$. (Looijenga
\cite{5}.)

Recently we have succeeded in extending the similar descriptions to
fourteen triangle singularities and six quadrilateral
singularities.

The following fourteen singularities are called triangle
singularities or exceptional singularities. (Arnold \cite{1}.)

$$\matrix
E_{\,12},&Z_{\,11},&Q_{\,10},\\
  &&\\
  E_{\,13},&Z_{\,12},&Q_{\,11},\\
  &&\\
  E_{\,14},&Z_{\,13},&Q_{\,12},\endmatrix
\quad\quad\quad
\matrix
W_{\,12},&S_{\,11},\\
  &\\
  W_{\,13},&S_{\,12},\endmatrix
\quad\quad\quad
U_{\,12}
$$

\vskip 10pt

For each $\Xi$ of the above fourteen classes the corresponding
Gabri\'{e}lov graph is defined. The type $(p,\,q,\,r)$ of the
corresponding Gabri\'{e}lov graph is as in the following list
(Gabri\'{e}lov \cite{4}):

$$\matrix
(2,\;3,\;7),&(2,\;4,\;5),&(3,\;3,\;4),\\
&&\\
(2,\;3,\;8),&(2,\;4,\;6),&(3,\;3,\;5),\\
&&\\
(2,\;3,\;9),&(2,\;4,\;7),&(3,\;3,\;6),\\
\endmatrix
\quad\quad\quad
\matrix
(2,\;5,\;5),&(3,\;4,\;4),\\
  &\\
(2,\;5,\;6),&(3,\;4,\;5),\endmatrix
\quad\quad\quad
(4,\;4,\;4)
$$

\vskip 10pt

\proclaim{Theorem 1}
{\rm (Urabe \cite{10}.)}\quad Let  $\Xi$  be one of the
above  fourteen classes of triangle singularities. The following two
conditions are equivalent.

\roster
\item"(A)"  $\Gamma\in PC (\Xi)$.

\item"(B)"  Either {\rm (B-1)} or {\rm (B-2)} holds.

\hbox{\qquad\quad\hfill\vbox{\hsize=11.5cm \roster
\item"(B-1)" $\Gamma$ can be made by an elementary
transformation or a tie transformation from a Dynkin subgraph of the
corresponding Gabri\'{e}lov graph to $\Xi$.\endroster}}
\newpage
\hbox{\qquad\quad\hfill\vbox{\hsize=11.5cm \roster
\item"(B-2)" $\Gamma$ is one of the following exceptions:

\quad For $\Xi=Z_{13}$, $A_7+A_4$.

\quad For $\Xi=S_{11}$, $2A_4+A_1$.

\quad For $\Xi=U_{12}$,
$2D_4+A_2$,$A_6+A_4$,$A_5+A_4+A_1$,$2A_4+A_1$.

\quad For the other eleven classes, no exceptions.
\endroster}}\endroster\endproclaim

A {\it tie transformation} in (B-1) is another operation by which we
can make a new Dynkin graph from a given Dynkin graph.

\proclaim{ Definition}\rm
(A tie transformation)\quad Assume that applying the following
procedure to a Dynkin graph  $\Gamma$, we have obtained the Dynkin
graph
${\overline\Gamma}$. Then we call the following procedure a tie
transformation of a Dynkin graph.

\roster
\item Add one vertex and a few edges to each component
of  $\Gamma$  and make it into the extended Dynkin graph of the
corresponding type. Moreover attach the corresponding coefficient of
the maximal root to each vertex.

\item Choose in an arbitrary manner subsets  $A$, $B$  of the set of
the vertices of the extended graph $\tilde \Gamma $ satisfying the
following conditions:\endroster

\hbox{\qquad\quad\hfill\vbox{\hsize=11.5cm \roster
\item"$\langle a\rangle$" $A\cap B=\emptyset $.

\item"$\langle b\rangle$" Let $V$  be the set of vertices of an
arbitrarily chosen component  $\tilde \Gamma '$
of  $\tilde \Gamma $.  Let  $\ell $
be the number of elements in  $V\cap A$ and
$n_{\,1},\;n_{\,2},\;...,\;n_{\,\ell }$ be the numbers attached
to $V\cap A$.  Furthermore let  $N$  be the sum of numbers attached
to
$V\cap B$.  (If  $V\cap B=\emptyset $, then  $N=0$ .)  Then the
greatest common divisor of the  $\ell +1$ numbers
$N,\;n_{\,1},\;n_{\,2},\;...,\;n_{\,\ell }$ is 1.
\endroster}}
\roster
\item[3] Erase all attached integers and remove
vertices belonging to  $A$  together with edges issuing from them.

\item[4] Draw another new vertex  $\bigcirc$  corresponding to a
root
$\alpha$  with $\alpha ^2=2$. Connect this new vertex  $\bigcirc$
and each vertex in  $B$  by an edge.  \endroster\endproclaim

\proclaim{Remark}\rm
Often the resulting graph  $\bar \Gamma $ after the above
procedure (1) -- (4) is not a Dynkin graph. We consider only the
cases where the resulting graph $\bar \Gamma $ is a Dynkin graph and
then we call the above procedure a tie transformation.  Under this
restriction the number $\#(B)$ of elements in the set  $B$  satisfies
$0\le \#(B)\le 3$.  $\ell =\#(V\cap A)\ge 1$.

Any book on Lie algebras contains the definition of integers in (1)
called {\it the coefficients of the maximal root}. (Bourbaki
\cite{2}.)
\endproclaim

\proclaim{Example}\rm
 We consider the case $\Xi=W_{13}$.
The Gabri\'{e}lov graph in this case is the following and it has a
Dynkin subgraph of type $E_8+A_2$:

\centerline{\picture 149mm by 22mm (W13 scaled 1000)}

First we apply a tie transformation to $E_8+A_2$.  In the second
step of the transformation we can choose subsets  $A$  and  $B$  as
follows:

\centerline{\picture 134mm by 56mm (E8A2 scaled 1000)}

For the component of type $E_8$, $\ell=1$, $n_1=4$, $N=1$, and
thus $G.C.D.(n_1,\;N)=1$.  For the component $A_2$, $\ell=1$,
$n_1=1$, $N=1$,and thus $G.C.D.(n_1,\;N)=1$. One sees that the
condition $\langle b\rangle$ is satisfied.  Adding a new vertex
$\theta$ in the fourth step, one gets a graph of type
$A_6+D_5$ as the result of the transformation.  By Theorem 1 one
can conclude $A_6+D_5\in PC(W_{13})$.

Second we apply an
elementary transformation to $E_8+A_2$.

\centerline{\picture 117mm by 44mm (E8A2-2 scaled 1000)}

As in the above figure we can make $E_6+2A_2$.  Thus $E_6+2A_2\in
PC(W_{13})$.
\endproclaim

Now, it is very strange that Theorem 1 has a few exceptions
in a few cases.  Perhaps this is because our theory has a missing
part.

\proclaim{Problem} \rm  Find the missing part of our theory and give
a simple characterization of the set $PC(\Xi)$ without exceptions.
\endproclaim

This problem may be very difficult, but I believe that there
exists a solution.

Now,  for the following nine singularities of the above fourteen
ones we can state another theorem:

$$\matrix
E_{\,12},&Z_{\,11},&Q_{\,10},\\
  &&\\
  E_{\,13},&Z_{\,12},&Q_{\,11},\\
  &&\\
  E_{\,14},&Z_{\,13},&Q_{\,12}.\endmatrix$$

\midspace{5mm}\pagebreak
\proclaim{Theorem 2}
{\rm (Urabe \cite{9}.)}\quad Let  $\Xi$  be one of the
above  nine classes of singularities.   The following two conditions
are equivalent.

\roster
\item"(A)"  $\Gamma\in PC (\Xi)$.

\item"(B)"  The Dynkin graph  $\Gamma$  has only  components of type
$A$,
$D$ or $E$, and can be made from the essential basic graph depending
on  $\Xi$ by a combination of two of elementary transformations and
tie transformations. \endroster

\vskip -5pt

$$\vbox{
\hbox{\qquad
The respective essential basic graph
}
\hbox{\it
corresponding to the above nine singularities
}
}$$
$$\matrix
E_{\,8},&E_{\,7},&E_{\,6}\\
  &&\\
  E_{\,8}+BC_{\,1},&E_{\,7}+BC_{\,1},&E_{\,6}+BC_{\,1}\\
  &&\\
  E_{\,8}+G_{\,2},&E_{\,7}+G_{\,2},&E_{\,6}+G_{\,2}
\endmatrix$$
\endproclaim

In the condition (B) four kinds of combinations --
``elementary" twice, ``tie" twice, ``elementary" after ``tie", and
``tie" after ``elementary" -- are all permitted.
Even when $\Xi=Z_{13}$ no exception appears in Theorem 2.
Note that Dynkin graphs of type $BC_1$ or $G_2$ appear
and the number of repetitions of transformations is not one but
two.

Here we give
some explanation on Dynkin graphs and root systems of type $BC$.  A
root system  $R$  is a finite subset of a Euclidean space satisfying
axioms on symmetry.  Usually we assume moreover the following axiom
(*) of the reduced condition:

\qquad\qquad\qquad\qquad\qquad\qquad
If  $\alpha \in R$, then
$2\alpha \notin
R$\qquad\qquad\hfill\hfill(*)

\noindent
Under these axioms we obtain irreducible root systems of type $A$,
$B$, $C$, $D$, $E$, $F$  and  $G$  as in any book on Lie algebras.
However, under the absence of the axiom (*) we have further a series
of irreducible root systems, which are called of type $BC_k$
$\left( {k=1,\;2,\;3,\;...} \right)$. (Bourbaki \cite{2}.)  It is
easy to generalize the concept of Dynkin graphs to root systems of
type
$BC$.  (Urabe \cite{8}.)  The Dynkin graph of type $BC_1$ is the
following: $\bigotimes$

We explain the meaning of this $BC_1$ graph.  Recall first the
meaning of Dynkin graphs.  Let  $R$  be an irreducible root system
and $\Delta \subset R$ be the root basis.  We can assume that the
longest root $\alpha\in R$ satisfies $\alpha ^2=2$ after normalizing
the inner product of the ambient Euclidean space.  The Dynkin graph
$\Gamma$ of  $R$  is the graph drawn by the following rules:  (1)
The vertices  of $\Gamma$ have one-to-one correspondence with the set
$\Delta$ (the root basis).  (2) Two vertices in $\Gamma$
corresponding to two elements $\alpha ,\;\beta \in \Delta $ are
connected by an edge in $\Gamma$ if and only if the inner product
$\left( {\alpha ,\;\beta } \right)\ne 0$.

If  $R$  is of type  $A$, $D$  or  $E$,  then  $R$  consists of only
roots $\alpha$ with $\alpha ^2=2$,  and every $\alpha \in \Delta $
satisfies $\alpha ^2=2$.  Therefore in these cases every vertex in
the Dynkin graph can be denoted by a small white circle $\bigcirc$.

If  $R$  is of type $BC_1$, then $\Delta$ consists of a unique root
$\delta$ with $\delta ^2=1/ 2$ and $R=\left\{ {\,-2\delta
,\;-\delta ,\;\delta ,\;2\delta \,} \right\}$. The vertex in the
Dynkin graph corresponding  to a root $\delta$ with $\delta ^2=1/
2$ is denoted by $\bigotimes$. The $BC_1$ graph is the graph
consisting of a unique vertex of this kind.  In this case the
maximal root $\eta $ is equal to $2\delta$, and thus the extended
Dynkin graph, i.e., the graph corresponding to $\Delta ^+=\Delta
\cup \left\{ {-\eta } \right\}$  is the
following:\enspace{\picture 34.3mm by 5.8mm (BC1 scaled 1000)}
\hskip -29pt  $\bigotimes$\hskip 19pt (The
edge is bold. The numbers are the coefficients of the maximal root.)

If  $R$  is of type $G_2$, then  $\Delta$ consists of two elements
$\alpha$  with $\alpha^2=2$ and $\gamma$ with $\gamma^2=2/3$.
We denote  the  vertex  corresponding  to a  root
$\gamma$   with  $\gamma^2=2/3$ by
{\picture 5.3mm by 5.7mm (G1 scaled 1000)}.
Our Dynkin graph of type $G_2$ is the following;
{\picture 24.7mm by 5.1mm (G2 scaled 1000)}
and our extended Dynkin graph of
type $G_2$ is the following (The numbers are the coefficients of the
maximal root.):

\centerline{\picture 45.9mm by 11mm (G2Ext scaled 1000)}

Note that as a result of an elementary or a tie transformation, a
graph  consisting of a unique vertex corresponding to $\gamma$ with
$\gamma^2=2/3$ can appear.  We call the graph
{\picture 5.3mm by 5.7mm (G1 scaled 1000)} the Dynkin graph of type
$G_1$. This corresponds to the root system
$R=\left\{ {\,-\gamma ,\;\gamma \,} \right\}$ with $\gamma^2=2/3$.
The extended Dynkin graph of type $G_1$ is the following:
\vbox to 15pt{}{\picture 35.3mm by 5.1mm (G1Ext scaled 1000)}.
(The edge is bold.
The numbers are the coefficients of the maximal root.)

We can explain why we do not use the standard
expression {\picture 25.8mm by 6.4mm (StandardG2 scaled 1000)}
of the $G_2$ graph. (Bourbaki \cite{2}.)  If we
use the standard expression, we cannot define the concept of the
$G_1$ graph.

Note that since we have assumed that the Dynkin graph $\Gamma$ in
Theorem 2 has only components of type  $A$, $D$  or  $E$,  any Dynkin
graph with a component of type  $G_2$, $G_1$ or $BC_1$ made by two
transformations has no meaning, and is to be thrown out.

In the next hierarchy of Arnold's list six kinds of quadrilateral
singularities $J_{3,0}$, $Z_{1,0}$, $Q_{2,0}$, $W_{1,0}$,
$S_{1,0}$, $U_{1,0}$ appear. Also for them we can show similar
theorems to Theorem 2. This is the theme of Urabe \cite{8}.
Though the theorems for them are also quite simple, we have to
introduce several additional concepts such as Dynkin graphs of type
$B$ or
$F$, obstruction components and so forth. Since they have been
treated in Urabe
\cite{8}, we omit the further explanation for them here.

\enddocument